\documentclass[12pt]{article}
\usepackage[utf8]{inputenc}
\usepackage{amsmath,amssymb,amsthm,lscape,longtable,enumerate}
\usepackage[english]{babel}
\usepackage{indentfirst}
\usepackage[dvips]{graphicx}

\renewcommand{\le}{\leqslant}
\renewcommand{\ge}{\geqslant}

\newtheorem{theorem}{Theorem}
\newtheorem{lemma}{Lemma}

\title{On the stability of some Erd\H{o}s--Ko--Rado type results\footnote{This work is done 
under the financial support of the grant 15-01-00350 of 
Russian Foundation for Basic Research.}}

\author{Pyaderkin M.M.\footnote{Moscow State University, Mechanics and Mathematics Faculty, Department of Probability Theory.}}
\date{}

\begin{document}

\maketitle

\section{Introduction}

One of the central results in extremal set theory is due to Erd\H{o}s, Ko, and Rado: they
studied \textit{intersecting} families $\mathcal{A}$ (see \cite{KoRado}). A family of sets $\mathcal{A}$ is said to be $t$-intersecting
if $|A \cap B| \ge t$ for any $A, B \in \mathcal{A}$. We are interested in families of sets of the same size: we write $[n]$
for the set $\{1,2,\ldots,n\}$ and $[n]^{(r)}$ for the family of all the subsets of $[n]$ of size $r$. 
The classical Erd\H{o}s--Ko--Rado theorem says that if $\mathcal{A} \subset [n]^{(r)}$ is a
$t$-intersecting family, then $|\mathcal{A}| \le \binom{n-t}{r-t}$ for all $r$, $t$ and sufficiently large $n$.
It is clear that this bound is tight: one can consider a family $\mathcal{A}$ consisting of subsets
which contain all the elements $1, 2, \ldots, t$. Such families are called \textit{trivial} ones.

There are many results on the structure of the non-trivial intersecting families (see \cite{Fra}, \cite{Alsw}). 
 One of the important results states that the size of any non-trivial $t$-intersecting family is much 
 smaller than the size of the trivial one. The exact statement is described in Theorem \ref{alsw} (see Subsection 2.5). 

Recently Bollob\'as, Narayanan, and Raigorodskii introduced a random setting of this problem, provided $ t = 1 $ (see \cite{Nar}). To do so, one first needs
to reformulate the problem in terms of graph theory. Namely, consider classical \textit{Kneser}'s graph $K(n,r)$: for
two natural numbers $ r, n $ such that $r \le n / 2$, its vertices are all the subsets of $[n]$ of size $r$, and two such vertices are adjacent
if the corresponding subsets are disjoint. It is clear that a $1$-intersecting family forms an independent set in
$K(n,r)$, and so we are interested in the size of the largest independent set in $K(n,r)$. We denote
this size by $\alpha(K(n,r))$. The Erd\H{o}s--Ko--Rado theorem states that $\alpha(K(n,r))=\binom{n-1}{r-1}$.

Now let us delete each edge of the graph $K(n,r)$ with some fixed probability $p$ independently of each other. Quite
surprisingly, the independence number of such random subgraph $K_p(n,r)$ of the graph $ K(n,r) $ is, with high probability, {\it the same} as the 
independence number of the initial graph.
This phenomenon is called the \textit{stability} of the independence number.
To be more precise, we formulate the result of Bollob\'as, Narayanan, and Raigorodskii, who established a threshold
probability for the stability property.

\begin{theorem}
Let $p_c(n,r)=\frac{(r+1) \log n - r \log r}{\binom{n-1}{r-1}}$ for $r=r(n)=o(n^{1/3})$, $r \ge 2$, and let us fix some $\varepsilon > 0$.
Then as $n \rightarrow \infty$

$$\mathbb{P}\left[ \alpha\left(K_p(n,r)\right) = \binom{n-1}{r-1} \right] \rightarrow \begin{cases}
1, & \text{if } p \le 1 - (1+\varepsilon)p_c(n,r), \\
0, & \text{if } p \ge 1 - (1-\varepsilon)p_c(n,r).
\end{cases}$$
\end{theorem}

One of the preceeding results comes from a bit different point of view. 
In \cite{We}, \cite{We1} Bogolyubskiy, Gusev, Pyaderkin,
and Raigorodskii studied random subgraphs $G_p(n,r,s)$ of $G(n,r,s)$, where $G(n,r,s)$ is the
graph whose vertex set is $[n]^{(r)}$ and two vertices $A, B \in [n]^{(r)}$ are adjacent if and only if
they intersect in exactly $s$ elements. These graphs are called \textit{distance graphs}, as 
they have a very nice geometric interpretation. Indeed, consider an $n$-dimensional hypercube $\{0,1\}^n$, and
construct a graph whose vertices form a subset of the set of vertices of this cube which have exactly $r$ non-zero
coordinates, and there is an edge between two vertices if their scalar product is exactly $s$. These graphs are well-known in combinatorial geometry and 
coding theory (see \cite{Rai3}--\cite{PA}). Of course, $ K(n,r) = G(n,r,0) $. 

In their paper \cite{We}, the authors present a series of theorems concerning the independence number
 of the graphs $G_p(n,r,s)$. In particular, they prove the following theorem.

\begin{theorem}
Let $ r, s $ be fixed. Then with high probability $$\alpha(G_{1/2}(n,r,s)) \ge (1+o(1)) \frac{\binom{n}{r} (r-s)}{\binom{r}{s}
\binom{n-r}{r-s}} \log_2 n.$$
\end{theorem}

This provides a lower bound for the independence number. On the other hand, there is
a work \cite{Pya1} by the author of the present paper where an upper bound is proven.

\begin{theorem}
\label{gnrs}
Let $r,s,\varepsilon > 0$ be fixed. There exists a constant
$\delta=\delta(\varepsilon, r,s)$ such that with high probability 
$$\alpha\left( G_{1/2}(n,r,s) \right) \le (1 + \varepsilon) 
\alpha\left(G(n,r,s)\right) + \delta \binom{n}{s} \log_2 n.$$
\end{theorem}

We see that the upper bound and the lower bound are of the same order: in case $r \le 2s+1$
the value of $\alpha(G(n,r,s))$ is $\Theta(n^s)$ (see \cite{FraFur}), and in the opposite case $r > 2s+1$ we know
that $\alpha(G(n,r,s))=\binom{n-s-1}{r-s-1}$ for any sufficiently large $n$ (see \cite{FraFur}).
Since in case $r > 2s+1$ the value $ \binom{n-s-1}{r-s-1} $ dominates asymptotically (as $ n \to \infty $) the value $ \binom{n}{s} \log_2 n $ 
(see Theorem 3), in this case, we get an {\it asymptotic} stability.
The main result of the present paper is much stronger, though less general.

\begin{theorem}
\label{main}
Let $r \ge 4$ be constant. Then with high probability $\alpha\left(G_{1/2}(n,r,1)\right)=\alpha\left(G(n,r,1)\right).$
\end{theorem} 

So, as in case of Kneser's graphs, we see a stability result for $p=\frac{1}{2}$: even if we delete half of the edges
of the graph, the size of a maximum independent set with high probability does not change. There are some clues that
the same behaviour occurs always when $r>2s+1$, not only in case $s=0$ or $s=1$. If this is true, we get a very impressive 
``jump'' from the stability in case $ r > 2s+1 $ to its complete absence in case $ r \le 2s+1 $. Indeed, it follows from Theorems 2 and 3
that for $ r \le 2s+1 $, with high probability $ \alpha(G_{1/2}(n,r,s)) = \Theta(\alpha(G(n,r,s)) \log_2 n) $. 

At the same time, we will see from the proof of Theorem 4 that 
of course one can replace the edge probability $ 1/2 $ by any other constant and even by a slowly decreasing function. However, to find 
a threshold probability as in Theorem 1 would be quite difficult.

The proof of Theorem \ref{main} is essentially based on the proof of Theorem \ref{gnrs}. However, to prove
a stability result
we need a more careful analysis of the structure of the independent sets in $G(n,r,1)$.  Actually,
we prove that any sufficiently large independent set in $G(n,r,1)$ has some predefined structure, and this
is a result of independent interest.

The family $S_{xy}$ of all $r$-element subsets of $[n]$ containing two fixed elements $x$ and $y$ $(x \neq y)$ is called the \textit{star centered
at elements $x$ and $y$ over the set $[n]$}. Let $d(A)$ 
denote the maximal size of a subset of $A$ which is a part of some star: $d(A)=\max_{x \neq y} |S_{xy} \cap A|$.
Our goal is to bound the size of $A$ given the value of $d(A)$.

\begin{theorem}
\label{r5}
Let $r \ge 5$ and $r=r(n)=o(n^{1/8})$. If $A$ is an independent set in $G(n,r,1)$ with 
$d(A) \le \binom{tn-2}{r-2}$ for some $t \in (0,1)$, then $|A| \le \frac{1}{t}\binom{tn-2}{r-2} + o\left(\binom{n-2}{r-2}\right)$, where
$\binom{a}{b}$ is defined as $\frac{a(a-1)\ldots(a-b+1)}{b!}$ for non-integer $a$.
\end{theorem} 

The right-hand side of the inequality in the theorem above consists of two terms, and actually it does not matter
for us which one is greater.
Anyway the theorem means that if $d(A)$ is not large, then $A$ itself is not large. Conversly, if $A$ is a rather large subset it implies that
$d(A)$ is large, so $A$ contains a large part of some star.

Note that the restriction $ r = o\left(n^{1/8}\right) $ is apparently not the best possible. With some additional technical work, it can 
certainly be weakened considerably. 

In some sense, we can say that the bound in Theorem 5 is tight for constant $r$. Suppose $\frac{1}{t}=k$ for some integer $k$ which also divides $n$,
and divide the ground set $[n]$ into $k$ equal parts $p_1=\{1,2,\ldots,tn\}$, $p_2=\{tn+1,tn+2,\ldots,2tn\}$,~\dots.
For each part $p_i$ consider a star $A_i$ over $p_i$ centered at some two elements, and let $A$ be the union of these stars.
Clearly, $|p_i|=tn$, so $d(A)=|A_i|=\binom{tn-2}{r-2}$ and $|A|=\sum_{i=1}^k |A_i|=\frac{1}{t}\binom{tn-2}{r-2}$, which is
asymptotically the same as the upper bound in the theorem above.

Somehow surprisingly, we need a bit more careful estimate in case $r=4$ (although the idea of the proof is the same).

\begin{theorem} 
\label{r4}
Let $r=4$. If $A$ is an independent set in $ G(n,4,1) $ with 
$d(A) \le \frac{(tn)^2}{2}$ for
some $t \in (0,1)$, then $|A| \le (1+o(1))\max\left(\frac{n^2}{8},\frac{tn^2}{2}\right)$.
\end{theorem} 

The unusual part in the upper bound, $\frac{n^2}{8}$, comes from another interesting example.
Let $n$ be even and let us split the ground set into $\frac{n}{2}$ disjoint pairs. One
can obtain a quadruple by merging two such pairs, and let $A$ be the set of all such quadruples.
It is easy to see that $|A|=\binom{n/2}{2}=(1+o(1))\frac{n^2}{8}$, and $d(A)=O(n)$, so the 
constructed set is rather large, despite the fact it does not contain a large part of any star.

We organize our paper in the following way:
in the next section we discuss the proofs of Theorems \ref{main}, \ref{r5} and \ref{r4}. 
Each proof uses some auxiliary facts (lemmas), which are proven
in the third section. 

\section{Proofs of Theorems \ref{main}, \ref{r5}, \ref{r4}}

\subsection{Proof of Theorem \ref{main}}

To prove Theorem \ref{main} we need to recall some definitions from the paper \cite{Pya1}.
Let us fix some positive real $\gamma$ and let $m=\left\lfloor\frac{6r\log_2n}{\gamma}\right\rfloor$ (note that $r$ is constant
in this theorem). We say that an ordered pair $(A_1, A_2)$, where $A_1$ and $A_2$ are two subsets of the set of vertices of $G(n,r,s)$,
\textit{forms a $\gamma$-dense construction}, if $|A_1|\ge m$, $|A_2| \ge m$
and each vertex from $A_2$ is adjacent with at least $\gamma |A_1|$ vertices from $A_1$. 
Note that, $A_1$ and $A_2$ can intersect each other or even coincide.

We say that $A$ \textit{contains a $\gamma$-dense construction}, if there is a pair of subsets $A_1 \subset A$, $A_2 \subset A$, 
such that $(A_1,A_2)$ forms a $\gamma$-dense construction.
We denote by $\mathcal{A}_\gamma$ the family of all such sets $A$, which contain a $\gamma$-dense construction.
Now we use Propositions 1 and 2 from \cite{Pya1} and we formulate them as the two theorems below.

\begin{theorem}
\label{prop1}
A subset of the set of vertices of $G(n,r,s)$, which contains a $\gamma$-dense construction, 
with high probability is not independent in $ G_{1/2}(n,r,s) $:
$$\mathbb{P} \left[ \exists A \in \mathcal{A}_\gamma: A \text{ is an independent set in } 
G_{1/2}(n,r,s) \right] \rightarrow 0.$$
\end{theorem}

\begin{theorem}
\label{prop2}
A subset of  the set of vertices of $G(n,r,s)$, which does not contain a  $\gamma$-dense construction, is
rather small, i.e. there exist two positive real functions $\mu_1(\gamma)$ and $\mu_2(\gamma)$ such that
$$\forall A \notin \mathcal{A}_\gamma: |A| \le \mu_1(\gamma) \binom{n}{s} \log_2n + \mu_2(\gamma) \alpha\left(A\right),$$
where $\mu_2(\gamma) \rightarrow 1$ as $\gamma \rightarrow 0$, 
and  $\alpha\left(A\right)$ stands for the size of a maximum independent subset of $A$ in $G(n,r,s)$.
\end{theorem}

We need to estimate the probability that in $G_{1/2}(n,r,1)$, $r \ge 4$, there is an independent set
of size at least $M+1$, where $M=\alpha\left(G(n,r,1)\right)=\binom{n-2}{r-2}$, provided $ n $ is large enough. To this end, we
need a lemma.

\begin{lemma}
\label{mainpart}
There exists $t_0 \in \left(\frac{1}{2}, 1\right)$ such that if $r = const \ge 4$, $|A|>M$ and $d(A) > \binom{t_0n-2}{r-2}$, 
then with high probability $A$ is not independent in $G_{1/2}(n,r,1)$.
\end{lemma}

We are going to prove this lemma in the third section (Subsection 3.5). 
Taking $t_0$ from Lemma \ref{mainpart} we fix some positive $\varepsilon \le \frac{1-t_0}{2t_0}$.
Recalling Theorem \ref{prop2} we see that there exists $\gamma_0$ such that $\mu_2(\gamma_0) < 1 + \varepsilon$.
Now we bound the probability that there is an independent set in $ G_{1/2}(n,r,s) $ of size strictly greater than $M$:


$$\mathbb{P}\left[ \exists A: |A|>M, A \text{ is an independent set in } G_{1/2}(n,r,1) \right] =$$
$$=\mathbb{P}\left[ \exists A \in \mathcal{A}_{\gamma_0}: \ldots \right]
+\mathbb{P}\left[ \exists A \notin \mathcal{A}_{\gamma_0}, d(A) > \binom{t_0n-2}{r-2}: \ldots \right]
+\mathbb{P}\left[ \exists A \notin \mathcal{A}_{\gamma_0}, d(A) \le \binom{t_0n-2}{r-2}: \ldots \right],$$
where the first term tends to zero due to Theorem \ref{prop1} and the second term tends to zero as well due to Lemma \ref{mainpart},
so we only need to estimate the third term, and we are going to prove
that this term is actually zero for sufficiently large values of $n$. If $d(A) \le \binom{t_0n-2}{r-2}$, then
applying directly Theorem \ref{r5} in case $r \ge 5$ or Theorem \ref{r4} in case $r=4$ 
we see that $\alpha(A) \le \frac{1}{t_0} \binom{t_0n-2}{r-2}+o\left(\binom{n-2}{r-2}\right)$.
Introducing $\delta=\mu_1(\gamma_0)$ we apply Theorem \ref{prop2}:

$$\forall A \notin \mathcal{A}_{\gamma_0}: |A| \le 
\mu_1(\gamma_0) \binom{n}{1} \log_2n + \mu_2(\gamma_0) \alpha(A) \le
\delta \binom{n}{1} \log_2n + (1+\varepsilon) \alpha(A) \le $$
 $$\le \delta n \log_2n + \left(1+\frac{1-t_0}{2t_0}\right) \frac{1}{t_0} \binom{t_0n-2}{r-2} + 
o\left(\binom{n-2}{r-2}\right) \le
\frac{1+t_0}{2t_0} \frac{1}{t_0} t_0^{r-2} \binom{n-2}{r-2} + o\left(\binom{n-2}{r-2}\right) \le$$
$$
\le \frac{1+t_0}{2}\binom{n-2}{r-2} + o\left(\binom{n-2}{r-2}\right) < M,$$
where the last inequality holds true for sufficiently large $n$.

Theorem 4 is proven. 

\subsection{Preliminaries for Theorems \ref{r5} and \ref{r4}}

As before, in this section and in what follows in most inequalities we assume that $n$ is sufficiently large.

In any independent set $A$ of any graph $ G(n,r,1) $ there is a maximum subset $A_0=\{v_1,v_2,\ldots,v_k\}$ such that none of its
vertices intersect each other. By $I_0$ we denote the union $v_1 \cup v_2 \cup \ldots \cup v_k$ of elements, which
are used in $A_0$.

Any vertex from $A \setminus A_0$ intersects at least one from $A_0$, otherwise $A_0$ is not the maximal one.
We denote by $A_i$, $i = 1, 2, \ldots$, the set of vertices which intersect exactly $i$ vertices from $A_0$:

$$A_i=\{v \in A \setminus A_0: \  |\{j: v \cap v_j \neq \emptyset\}| = i\}.$$

Note that $i$ can not be greater than $k$, as there are exactly $k$ vertices in $A_0$. It can not be greater than
$r / 2$, as if vertex $v$ intersects $u \in A_0$ then it must intersect it by at least two elements. Let us
denote the maximal index $i$, for which $A_i$ is not empty, by $q$.
We also say that element $x \in [n] \setminus I_0$ is \textit{connected} to $v_i \in A_0$ if there are at least
$\omega=\omega(r,n)$ vertices from $A_1$ containing element $x$ and intersecting $v_i$, where

$$\omega(r,n)=\begin{cases}
1 & \text{for }r = 4, \\
r^5 \binom{n}{r-5} & \text{for }r \ge 5.
\end{cases}$$                    

We say that two different elements $x$ and $y$ are \textit{joint} if

$$\forall v \in A_1 \ \  |v \cap \{x, y\}| \neq 1,$$
so if $v$ contains $x$ it must contain $y$ as well and vice versa. We denote the set of
all elements which are joint with some other elements from $[n] \setminus I_0$ by $P$. Obviously, $P$ can be split into disjoint equivalence classes
$P=p_1 \sqcup p_2 \sqcup \ldots \sqcup p_l$, where in each class elements are pairwise joint.
Now we are ready to formulate a simple lemma, which is very useful to describe the structure of independent sets.

\begin{lemma}
\label{onevertex}
Element $x \in [n] \setminus (I_0 \cup P)$ can not be connected to two different vertices.
\end{lemma} 

Lemma \ref{onevertex} will be proven in the third section (Subsection 3.1) and it 
means that we can split some of the remaining elements into disjoint sets $J_i$, each containing the elements which are
connected to the vertex $v_i$:

$$J_i = \{x \in [n] \setminus (I_0 \cup P): |\{v \in A_1: x \in v, v \cap v_i \neq \emptyset\}| \ge \omega \}.$$

Note that there actually might be some elements which are not connected to any vertex $v_i$, and thus they do not belong
to any of the sets $J_i$.

It is almost clear that sets $A_i$ with $i \ge 2$ are not important at all for our studies (see Lemma 4 in case of Theorem 5 and Lemma 5 in case of
Theorem 6 in the two corresponding subsections below), and using the constraint for $d(A)$ (given in the statements of both Theorems 5 and 6)
we get another lemma to bound the size of $A_1$.

\begin{lemma}
\label{a1}
If $ r = o\left(n^{1/8}\right) $ and $d(A) \le \binom{tn-2}{r-2}$, 
then $|A_1| \le o\left(\binom{n-2}{r-2}\right) + 2kl + \frac{|I_0|+\sum_{i=1}^k |J_i|}{tn - 2}{\binom{tn-2}{r-2}}$.
\end{lemma}

A proof of this lemma is in Subsection 3.2. 

\subsection{Proof of Theorem \ref{r5}}

For $r \ge 5$, we first justify that we do not need to consider $A_i$ for $i \ge 2$.

\begin{lemma}
\label{ai}
For $r \ge 5$, $ r = o\left(n^{1/8}\right) $, we have $\sum\limits_{i=2}^q |A_i| = o\left(\binom{n-2}{r-2}\right)$.
\end{lemma} 

A proof is given in Subsection 3.3, and now we are ready to prove Theorem \ref{r5}.
If $A$ is an independent set in the graph $G(n,r,1)$ with $r \ge 5$, we can
apply Lemma \ref{ai} and bound the size of $A$:

$$|A|=\sum_{i=0}^q |A_i|=|A_0|+|A_1|+\sum_{i=2}^q |A_i| \le n + |A_1| + o\left(\binom{n-2}{r-2}\right)=|A_1| + o\left(\binom{n-2}{r-2}\right).$$

Now given $d(A) \le \binom{tn-2}{r-2}$ we use Lemma \ref{a1} to complete the proof:

$$|A| \le |A_1| + o\left(\binom{n-2}{r-2}\right) \le 2kl + \frac{|I_0|+\sum_{i=1}^k |J_i|}{tn - 2}\binom{tn-2}{r-2} + o\left(\binom{n-2}{r-2}\right) \le $$
$$ \le \frac{n}{tn-2}\binom{tn-2}{r-2} + o\left(\binom{n-2}{r-2}\right) \le \frac{1}{t}\binom{tn-2}{r-2} + o\left(\binom{n-2}{r-2}\right).$$

\subsection{Proof of Theorem \ref{r4}}

If a vertex $v$ intersects vertex $v_i$, then they share at least two elements, and this
implies, in case $r=4$, that the maximal index $q$, for which $A_i$ is not empty, does not exceed $2$. It means that
we need to deal only with $A_1$ and $A_2$, and since the size of the first one
is bounded by Lemma \ref{a1}, we now estimate the size of the second one.

\begin{lemma}
\label{a2}
If $r=4$, then $|A_2| \le 2k^2 + o(n^2)$.
\end{lemma}

A proof is given in Subsection 3.4. All parts of $A$ are bounded now:

$$|A|=|A_0|+|A_1|+|A_2| \le k + 2k^2 + 2kl + \frac{|I_0|+\sum_{i=1}^k |J_i|}{tn-2}{\binom{tn-2}{2}} + o(n^2).$$

Since all $J_i$ are pairwise disjoint, we have $\sum_{i=1}^k |J_i| \le n-|I_0|-|P| \le n-4k-2l$ and thus

$$|A| \le k+2k^2 + 2kl + \frac{n-4k-2l}{tn-2}\binom{tn-2}{2} + o(n^2) \le \frac{1}{2} \left( (2k+l)^2 + tn(n-4k-2l) \right) + o(n^2).$$

Now we see that the expression above depends only on $\Theta=\frac{4k+2l}{n} \in [0,1]$, and we rewrite the bound 
$$|A| \le \frac{n^2}{2} \left( \frac{1}{4} \Theta^2 + t(1-\Theta) \right) + o(n^2) \le \frac{n^2}{2} \max\left(\frac{1}{4}, t\right) + o(n^2) =
\max\left(\frac{n^2}{8}, \frac{tn^2}{2}\right) + o(n^2),$$
which completes the proof of Theorem \ref{r4}.

\section{Proofs of lemmas}

\subsection{Proof of Lemma \ref{onevertex}}

For any element $x \not \in I_0 \cup P$, we need to prove that there can exist at most one vertex $v_i$ connected to it.
Consider the set of all vertices from $A_1$ containing $x$. Each of them
interesects exactly one of the vertices from $A_0$, so we consider the sets $B_i$:

$$B_i=\{v \in A_1: x \in v, v \cap v_i \neq \emptyset\}.$$

For each vertex $v \in B_i$, we construct a set $f(v)$, which is obtained
by removing from $v$ two elements from $v_i$ and element $x$:

$$f(v)=v \setminus \{x, y, z\}, \text{where } y, z \in v \cap v_i.$$

If there are more than two elements in $v$ which belong to $v_i$, we remove any two of them.
Note that, $|f(v)|=r-3$.

\begin{lemma}
\label{difnon1}
If $B_i \neq \emptyset$ and $B_j \neq \emptyset$, then for any $u_1 \in B_i, u_2 \in B_j$, we have
$f(u_1) \cap f(u_2) \neq \emptyset$.
\end{lemma}

\noindent{\textbf{Proof of Lemma \ref{difnon1}.}}
Any two vertices $u_1 \in B_i$ and $u_2 \in B_j$ intersect in element $x$. Therefore, they must share at least one more
element, and this element can not belong to $v_i$: otherwise $u_2$ intersects two different vertices 
$v_i$ and $v_j$ from $A_0$. Obviously, it can not be an element from $v_j$ either. Thus, this element belongs to
$f(u_1) \cap f(u_2)$ and this completes the proof.
\vskip+0.2cm

In case $r=4$ for each vertex $u$ we know that $f(u)$ consists of one element, 
and it follows from Lemma \ref{difnon1} that
actually $f(u_1)=f(u_2)$ for any two vertices from different subsets $B_i$ and $B_j$.
If an element $x \in [n] \setminus (I_0 \cup P)$ is connected to two different vertices, then 
at least two among the sets $B_i$ are not empty, and this implies that $f(u_1)=f(u_2)$
for any two vertices (including vertices belonging to the same set $B_i$).
In this case there exists $y$ such that each vertex containing $x$ contains $y$ as well.
Applying the same argument to $y$ we derive that elements $x$ and $y$ are joint, but it contradicts
the assumption that $x \notin P$.

\vskip+0.2cm

In case $r \ge 5$ we need a more careful analysis.
Quite surprisingly, if we are given a non-empty set $ B_i $, then the existence of another large set $B_j$ 
yields some properties of the set $B_i$.

\begin{lemma}
\label{difnon2} If $B_i \neq \emptyset$ and $|B_j| \ge \omega$, then for any $u_1, u_2 \in B_i$, we have $f(u_1) \cap f(u_2) \neq \emptyset$.
\end{lemma}

\noindent{\textbf{Proof of Lemma \ref{difnon2}.}}
Suppose there exist two vertices $u_1, u_2 \in B_i$ such that $f(u_1) \cap f(u_2) = \emptyset$.
Each vertex $u$ from $B_j$ contains element $x$, and thus it intersects
$u_1$ in at least one element. Therefore, $u$ intersects $u_1$ in at least two elements, so there exists $y_1 \in u_1 \cap u \setminus \{x\}$.
It is easy to see that $y_1 \notin v_i$, otherwise $u$ intersects both $v_i$ and $v_j$, which contradicts that $u \in A_1$.
From that we conclude that $y_1 \in f(u_1) \cap u$. 
In the same manner one can show that $y_2 \in f(u_2) \cap u$, and $y_1 \neq y_2$, according to the initial 
assumption that $f(u_1) \cap f(u_2) = \emptyset$.

It follows that each $u \in B_j$ contains at least one element from $f(u_1)$ and at least one element from $f(u_2)$. 
It must contain element $x$ as well and intersects $v_j$ by at least two elements.
This implies that the size of $B_j$ is bounded by

$$|f(u_1)| |f(u_2)| \binom{r}{2} \binom{n}{r-5} \le r^4 \binom{n}{r-5} < \omega,$$
and it means that $|B_j| < \omega$.
This contradiction completes the proof of Lemma \ref{difnon2}.

\vskip+0.2cm

Suppose that element $x \in [n] \setminus (I_0 \cup P)$ is connected to two different vertices and consider $C_i = \{f(v): v \in B_i \}$. It means
that at least two of the sets $B_i$ have size not less than $\omega$.
From Lemmas \ref{difnon1} and \ref{difnon2} one can derive that any two subsets from $C=\bigcup_{i=1}^k C_i$ intersect each other.
We know that each vertex from $B_i$ can be obtained by adding two elements from $v_i$ to some element of $C_i$, so
$|B_i| \le \binom{r}{2} |C_i|$. Since the size of at least one among the sets $B_i$ is greater than $\omega$, 
we can bound the size of $C$:

$$|C| \ge |C_i| \ge \frac{|B_i|}{\binom{r}{2}} \ge \frac{\omega}{\binom{r}{2}} > 2 r^2 \binom{n}{r-5} = \omega'.$$

Now we are back to the Erd\H{o}s--Ko--Rado case: we have a 1-intersecting family $C$. We would like to prove that
our family is actually a part of some star. To do so, we use a theorem by Ahlswede and Khachatrian, which
they call ``The Complete Nontrivial-Intersection Theorem
for Systems of Finite Sets'' (see \cite{Alsw}).

\begin{theorem}
\label{alsw}
Let $A$ be a $t$-intersecting family of $k$-element subsets of $[n]$, and suppose this family is not a trivial one (is not a part
of some star). For any $ a, b \in {\mathbb N} $ such that $ a \le b $, let $ [a,b] = \{a, a+1, \dots, b\} $.
If $n > (t+1)(k-t+1)$, then $|A| \le \max\{|\mathcal{V}_1|,|\mathcal{V}_2|\}$, where

$$\mathcal{V}_1=\left\{v \in [n]^{(k)}: ~ |[1,t+2] \cap v| \ge t+1\right\},$$
$$\mathcal{V}_2=\left\{v \in [n]^{(k)}: ~ [1,t] \subset v, v \cap [t+1,k+1] \neq \emptyset \right\} \cup
\{[1,k+1] \setminus \{i\}: i \in [1,t]\}.$$
\end{theorem}

We are going to apply this theorem for $k=r-3$ and $t=1$: one can easily see that $$|\mathcal{V}_1| \le 3 \binom{n-3}{r-5} + \binom{n-3}{r-6} < \omega',$$
and $$|\mathcal{V}_2| \le \sum_{i=1}^{r-3} \binom{r-3}{i} \binom{n-(r-3)-1}{(r-3)-i-1} + 1 < (r-3) \binom{r-3}{1} \binom{n-(r-3)-1}{(r-3)-1-1} < \omega'.$$ 
Here and in what follows, we use the fact that, under the restrictions of the lemmas, the first summand is the maximum one, so that the sum is bounded by 
it times the number of summands. We do not check this explicitely, for this is a standard computation. 

As $|C| \ge \omega'$ 
it follows that $C$ is a part of a star, so 

$$\exists y \notin I_0: \forall v \in A_1, x \in v \implies y \in v.$$
Thus, all vertices containing $x$ contain $y$ as well.
As element $x$ is connected to at least two different vertices, so does element $y$. 
We apply the same argument to $y$ and derive that

$$\exists z \notin I_0: \forall v \in A_1, y \in v \implies z \in v.$$

If $z \neq x$ then each vertex containing $x$ contains two elements $y$ and $z$ as well.
It follows that we can bound the size of the set $B_i$: 
each vertex from this set contains elements $x$, $y$, $z$ and intersects $v_i$ by at least two elements. Thus,
$|B_i| \le \binom{r}{2} \binom{n}{r-5} < \omega$, which contradicts the
assumption that $x$ is connected to $v_i$.
So, the inequality is false and $z=x$ holds, and thus elements $x$ and $y$ are joint. This implies that $x \in P$ and
this contradiction completes the proof of Lemma \ref{onevertex}.

\subsection{Proof of Lemma \ref{a1}}

To bound the size of the set $A_1$, we need to consider a set $T \subset A_1$ of vertices, 
which contain some element $x \in [n] \setminus \left(I_0 \cup P\right)$
and intersect some vertex $v_i \in A_0$ not connected to element $x$:

$$T=\{v \in A_1: \exists x \notin (I_0 \cup P), \exists v_i \in A_0, x \text{ is not connected to } v_i, x \in v, v \cap v_i \neq \emptyset \}.$$
Notice that $T$ is an empty set in case $r=4$ because $\omega=1$. 

We now bound the size of the set $A_1$ by estimating the sizes of the set $T$ and of the remaining part:

\begin{equation}
|A_1|=|T|+|A_1 \setminus T|=|T|+|\{v \in A_1 \setminus T: v \cap P \neq \emptyset\}|+|\{v \in A_1 \setminus T: v \cap P = \emptyset\}|.
\end{equation}

\paragraph{An estimate of the first term in expression (1).}
In case $r=4$ set $T$ is empty, so we now assume that $r \ge 5$.
By definition of $T$, for each of its vertices, there is a vertex $v_i \in A_0$
which intersects it and an element $x$ which is not connected to this vertex. Thus,

$$|T| \le \sum_{x, v_i: \ x \text{ is not connected to } v_i} |\{v \in T: x \in v, v \cap v_i \neq \emptyset\}| \le 
\sum_{x, v_i: \ x \text{ is not connected to } v_i} \omega,$$
as element $x$ is not connected to the vertex $v_i$. Summing over all possible pairs $(x,v_i)$
and taking into account the condition $ r = o\left(n^{1/8}\right) $ we get the inequality

$$|T| \le n k \omega = o\left(\binom{n-2}{r-2}\right).$$

\paragraph{An estimate of the second term in expression (1).}
The second term estimates the number of vertices which intersect $P$. Each vertex intersecting $P$ can
intersect exactly one equivalence class from $P$ or more than one class. We first deal with the second case.

In cases $r=4$ and $ r = 5 $, there are no such vertices: such vertex must contain at least $6$ elements.
Indeed, if a vertex crosses some class $p_i$, as all the elements in the same class are joint, it means
that this vertex must contain the entire class. In each class there are at least two elements, and
each vertex must intersect one vertex $v_i$ from $A_0$ in at least two elements. 

Let $ r \ge 6 $ and take some $ v \in A_1 \setminus T $ such that $ v $ intersects $ P $ by at least two classes.
Then the total number of such vertices $ v $ for fixed $i$ is at most $\binom{r}{2} \binom{l}{2} \binom{n}{r-6}$:
we choose two classes from $P$, we choose two elements from the appropriate $v_i$ and we choose the remaining
elements, at most $r-6$, from $[n]$.
Summing over all possible $i$ and taking into account the condition $ r = o\left(n^{1/6}\right) $ 
(which is even weaker than the condition of the lemma) we get

$$|A_0| \binom{r}{2} \binom{l}{2} \binom{n}{r-6} \le nr^2n^2 \binom{n}{r-6} = o\left(\binom{n-2}{r-2}\right).$$

The remaining vertices intersect $P$ in exactly one class. In case $r \ge 5$ a very simple bound suffices.
Let us fix some class $p_j$ and consider
all vertices from $A_1 \setminus T$ intersecting this class. If a vertex intersects a class, then it must contain
it entirely. Also each vertex  $v \in A_1 \setminus T$ must intersect at least one vertex $v_i \in A_0$. 
If $v$ intersects some $v_i$, then it can contain, besides elements from $I_0 \cup P$, only elements from $J_i$.
Indeed, if $y \notin I_0 \cup P$ and $y \in v$, then as $v \notin T$
element $y$ is connected to $v_i$, so $y \in J_i$. It follows that those vertices, intersecting both $p_j$ and $v_i$, contain only
elements from $v_i \cup p_j \cup J_i$. As these vertices contain the entire class $p_j$ and at least two
elements from $v_i$, there are at most $\binom{r}{2} \binom{|v_i \cup J_i|}{r-4}$ of them.
A direct summation gives the following bound on the number of vertices intersecting $p_j$:

$$\sum_{i=1}^k \binom{r}{2} \binom{|v_i \cup J_i|}{r-4} \le r^2 \binom{\sum_{i=1}^k (|J_i| + r)}{r-4} \le r^2 \binom{n}{r-4}.$$

Summing over all the classes $p_j$ and taking into account the condition $ r = o\left(n^{1/4}\right) $ (which is once again 
weaker than the condition of the lemma) we bound the required number of vertices by $l r^2 \binom{n}{r-4}=o\left(\binom{n-2}{r-2}\right)$. 

The only thing left is dealing with $r=4$. Consider some vertex $v_i$.
We say that a pair $(x,y)$ of two elements $ x,y \in v_i$ is \textit{significant,} if

$$\exists j_1(x,y), j_2(x,y): ~ j_1(x,y) \neq j_2(x,y), \{x, y\} \cup p_{j_1(x,y)} \in A_1, \{x, y\} \cup p_{j_2(x,y)} \in A_1.$$

Why do we use such a complicated definition? 
First, two significant pairs can not intersect in one element. Indeed, suppose there are two 
significant pairs $(x,y)$ and $(x,z)$. These pairs intersect already in one element,
and classes $p_a$ and $p_b$ are disjoint when $a \neq b$, so $j_1(x,y)=j_1(x,z)$ and $j_1(x,y)=j_2(x,z)$, 
and thus $j_1(x,z)=j_2(x,z)$, which contradicts the definition of a significant pair.
Second, if a pair $(x,y)$ is not significant it means that there exists at most one class $p_j$
such that $\{x,y\} \cup p_j \in A_1$.

Since two significant pairs can not intersect each other, in each $v_i$ there are at most two significant pairs. Now the number of 
vertices which contain at least one significant pair and one class from $P$ is bounded by $2kl$.
The number of vertices which do not contain a significant pair is bounded by $k \binom{r}{2}=O(n)$, because
in each vertex $v_i$ there are $\binom{r}{2}$ pairs. 

\paragraph{An estimate of the third term in expression (1).}
Now consider the set
$B_i$ of all vertices from $A_1 \setminus T$, intersecting $v_i$ and not containing joint elements:

$$B_i=\{v \in A_1 \setminus T: v \cap v_i \neq \emptyset, v \cap P = \emptyset \}.$$

Each vertex from $B_i$ consists of elements from $v_i \cup J_i$ and does not intersect any of $v_j$, where $j \neq i$.
If there are two disjoint vertices $u_1$ and $u_2$ from $B_i$, then $A_0$ is not the maximal one:
we can discard $v_i$ from $A_0$ and add $u_1$ and $u_2$ instead.
This implies that in $B_i$ any two vertices intersect each other. They can not intersect in one element,
so they intersect in at least two. Now we are back to the Erd\H{o}s--Ko--Rado problem: we have a 2-intersecting
family. We use the result of Ahlswede and Khachatrian (Theorem \ref{alsw}) again, but now for $t=2$ and $k=r$.
Notice that we consider only vertices which consist of elements from $v_i \cup J_i$, and thus
in the theorem we are going to replace $n$ with $|v_i \cup J_i|=r+|J_i|=\tilde{n}.$
In this case we have $$|\mathcal{V}_1| \le 4 \binom{\tilde{n}-4}{r-3} + \binom{\tilde{n}-4}{r-4},$$ and $$|\mathcal{V}_2| \le 
2+\sum_{i=1}^{r-2} \binom{r-1}{i} \binom{\tilde{n}-r-1}{r-i-2}.$$

It is easily checked that in case $ \tilde{n} \ge r^2 $, the last sum is not greater than
$ r^2 \binom{\tilde{n}-r-1}{r-3} < r^2 \binom{\tilde{n}}{r-3} $. Thus, 
either $\tilde{n} < r^2$, or any non-trivial intersecting family has at most $r^2 \binom{\tilde{n}}{r-3}$ vertices. 
Before summing over all $B_i$, we need to distinguish 
three different cases.

\begin{enumerate}
\item $|J_i| + r < r^2$. In this case we bound the size of $B_i$ by $\binom{|J_i|+r}{r}$, and summing
over $i$ yields
$$\sum_{i} \binom{|J_i|+r}{r} \le \sum_{i} \binom{|J_i|+r}{r-3} \left(r^2\right)^3 \le \binom{\sum_{i} \left(|J_i|+r\right)}{r-3} (r^2)^3 
\le r^6 \binom{n}{r-3} = o\left(\binom{n-2}{r-2}\right).$$

\item $|J_i| + r \ge r^2$, but $B_i$ is not a trivial family. Here summing over $i$ yields at most
$r^2 \binom{n}{r-3}=o\left(\binom{n-2}{r-2}\right)$.

\item $B_i$ is a trivial family. This implies that on the one hand $|B_i| \le \binom{|J_i|+r-2}{r-2}$, but
on the other hand by the initial assumption we have $d(A) \le \binom{tn-2}{r-2}$, so actually

$$|B_i| \le \min\left(\binom{|J_i|+r-2}{r-2}, \binom{tn-2}{r-2}\right) = \binom{\min\left(|J_i|+r-2,tn-2\right)}{r-2}$$

for sufficiently large $n$. Summing over $i$ yields at most

$$
\sum_{i=1}^k \binom{\min\left(|J_i|+r-2,tn-2\right)}{r-2} \le \sum_{i=1}^k \binom{tn-2}{r-2} \frac{|J_i|+r-2}{tn-2}
\le \frac{kr + \sum_{i=1}^k |J_i|}{tn-2} \binom{tn-2}{r-2}.
$$

\end{enumerate} 

Lemma \ref{a1} is proven.

\subsection{Proof of Lemma \ref{ai}}

Recall that we need to estimate the size of $A_i$ for $i \ge 2$ that is to bound the number of
vertices which intersect at least two vertices from $A_0$.

Consider the set $A_i$. By definition each vertex $v \in A_i$ intersects
exactly $i$ vertices from $A_0$. Also, if $v$ intersects $v_j$, 
then they share at least two elements.
From that we conclude that the size of $A_i$ is not greater than
$\binom{|A_0|}{i} \binom{r}{2}^i \binom{n-2i}{r-2i}$. Summing over $i \ge 3$ yields

$$\sum_{i=3}^q \binom{k}{i} \binom{r}{2}^i \binom{n-2i}{r-2i} \le \sum_{i=3}^q \binom{n}{i} \binom{r}{2}^i \binom{n-2i}{r-2i} \le
r \binom{n}{3} \binom{r}{2}^3 \binom{n-6}{r-6} = o\left(\binom{n-2}{r-2}\right),$$
where we essentially use the condition $ r = o\left(n^{1/8}\right) $.

So the only remaining set is $A_2$. First, we consider $v \in A_2$ which do not contain
any elements from $[n] \setminus I_0$. These vertices are contained in the union of
two vertices from $A_0$, and thus there are at most $\binom{k}{2} \binom{2r}{r}$ of them.
We use the condition $ r = o\left(n^{1/8}\right) $ again to claim that the last estimate is $o\left(\binom{n-2}{r-2}\right)$:

$$\binom{k}{2} \binom{2r}{r} \le n^2 (2r)^5 (2r)^{r-5} \le \frac{n^3}{n} (2r)^5 \left(\frac{n}{r}\right)^{r-5} \le 
\frac{32r^8}{n} \left(\frac{n}{r}\right)^{r-2} \le o(1) \left(\frac{n-2}{r-2}\right)^{r-2} \le o(1) \binom{n-2}{r-2}.$$

Now let us fix an element $x \notin I_0$ and $v \in A_2$ such that $x \in v$.
Since $v \in A_2$, it follows that
 $v \cap v_i \neq \emptyset$ and $v \cap v_j \neq \emptyset$ for some $i \neq j$. Any other
 vertex $u \in A_2$ such that $x \in u$ must intersect $v$ in at least one more element $y \neq x$. Again, we disinguish two different cases.

\begin{enumerate}
\item $y \notin v_i$, $y \notin v_j$. This implies that $y \in v \setminus (v_i \cup v_j)$, and $y$ can not
belong to any vertex from $A_0$: otherwise $v \notin A_2$, since it intersects three or more vertices from $A_0$.
We now need to estimate the number of vertices from $A_2$ containing both $x$ and $y$.
Each such vertex must intersect two vertices from $A_0$, so the total number of such vertices is
less than $\binom{k}{2} \binom{r}{2}^2 \binom{n}{r-6} \le r^4 \binom{n}{r-4}$. As there are at most $r$ different elements $y$, it means that
the number of vertices from $A_2$ containing element $x$ and some element $y$, is bounded by $r^5 \binom{n}{r-4}$.

\item $y \in v_i$ (one can deal with the case $y \in v_j$ similarly). As $y \in v_i$, vertex $u$ already intersects $v_i$.
It must intersect $v_i$ in one more element, and it must intersect one more vertex from $A_0$. 
Obviously, it must contain element $x$ as well. The total number
of vertices satisfying such conditions is less than $n \binom{r}{2}^2 \binom{n}{r-5} \le r^5 \binom{n}{r-4}$.
\end{enumerate}

In both cases we see that the total numbers of vertices $v \in A_2$ containing element $x$ is not greater than $r^5 \binom{n}{r-4}$. Summing
over all possible elements $x$ we derive that the total size of $A_2$ is less than $o\left(\binom{n-2}{r-2}\right)$, and this
completes the proof of Lemma \ref{ai}.

\subsection{Proof of Lemma \ref{a2}}

In this lemma we need to bound the size of the set $A_2$ in case $r=4$.
We note that in case $r=4$ each vertex from $A_2$ consists of two pairs,
where one pair belongs to one vertex from $A_0$ and the other pair belongs to another vertex from $A_0$.

We organize the proof in the same way as it was done before: we say
that a pair $(x,y)$ of elements $x,y \in v_i$ is \textit{significant,} if

$$\exists u_1, u_2 \in A_2, j_1(x,y), j_2(x,y): \{x,y\} \subset u_1, u_1 \cap v_{j_1} \neq \emptyset, 
\{x,y\} \subset u_2, u_2 \cap v_{j_2} \neq \emptyset.$$

The key observation is that two significant pairs can not intersect in one element. Indeed, if
for example, $(x,y)$ and $(x,z)$ are significant, then they share at least one element,
and as vertices $v_a$ and $v_b$ are disjoint for $a \neq b$, so $j_1(x,y)=j_1(x,z)$ and $j_1(x,y)=j_2(x,z)$, 
and this implies that $j_1(x,z)=j_2(x,z)$,
which contradicts the definition of a significant pair.

Using the observation from the paragraph above,
we conclude that for each $i$ there are at most two significant pairs. That means 
that the total number of vertices from $A_2$, which consist of significant pairs, is less than
$\binom{2k}{2} \le 2k^2$. On the other hand, the total number of vertices which do not contain a significant pair, is less than
$k \binom{r}{2} \binom{r}{2}=O(n)$, and this completes the proof.

\subsection{Proof of Lemma \ref{mainpart}}
Let $t$ be a fixed positive real number and $d(A) > \binom{tn-2}{r-2} \ge t^{r-2}M+o(M)$. In this case, there
are two elements $x$ and $y$ such that
$|A \cap S_{xy}|>t^{r-2}M+o(M)$. As we know that $|A|>M$, 
there exists a vertex $v \in A$, which is not contained in $S_{xy}$. We can say that $x \notin v$.
Vertex $v$ is adjacent in $G(n,r,1)$ to every vertex of $S_{xy}$, excluding vertices
containing at least two elements from $v$ and element $x$. 
The total number of vertices from $S_{xy}$ which are not adjacent to $v$ is less than $\binom{|v|}{2} \binom{n}{r-3}=o(M)$,
so there are at least $t^{r-2}M + o(M)$ vertices adjacent to $v$, and the probablity that in $G_{1/2}(n,r,1)$ all these
edges are absent is not greater than $2^{-t^{r-2}M + o(M)}$.

Now we bound the number of ways to choose a star $S_{xy}$, its subset of size $t^{r-2}M+o(M)$, 
and a vertex $v \notin S_{xy}$, so that the probability is bounded by

$$\binom{n}{r} \binom{n}{2} \binom{M}{t^{r-2}M+o(M)} 2^{-t^{r-2}M + o(M)} = \binom{n}{r} \binom{n}{2} \binom{M}{(1-t^{r-2})M+o(M)} 2^{-t^{r-2}M + o(M)}$$

Using quite standard inequality $\binom{n}{k} \le \left(\frac{en}{k}\right)^k$ we continue the estimation:

$$n^r n^2 \left(\frac{eM}{(1-t^{r-2})M+o(M)}\right)^{(1-t^{r-2})M+o(M)} 2^{-t^{r-2}M + o(M)} = 2^{(1+o(1))\left(-(1-t^{r-2})M \log_2\left((1-t^{r-2})/e\right) - t^{r-2}M\right)}=$$
$$=2^{(1+o(1))M\left(-(1-t^{r-2})\log_2\left((1-t^{r-2})/e\right)-t^{r-2}\right)}.$$

As $\lim\limits_{x\rightarrow1}-(1-x) \log_2\left((1-x)/e\right) - x=-1$, so there exists $t_0 \in (\frac{1}{2},1)$, such that
the required probability is less than $2^{-(1+o(1))\frac{1}{2}M} \rightarrow 0,$
and this completes the proof of Lemma \ref{mainpart}.

\end{document}